\documentclass[10pt,oneside,a4paper,draft]{article} 	

\usepackage[fleqn]{amsmath} 						
\usepackage{amsthm,amsfonts,latexsym,amssymb,amscd} 	
\usepackage[mathscr]{eucal}   						
\usepackage[all,2cell]{xy} \UseAllTwocells 		
\usepackage{txfonts}

\usepackage[draft=false]{hyperref} 								

\renewcommand{\H}{\mathcal{H}} 	
 
\newcommand{\CC}{{\mathbb{C}}}

\newcommand{\NN}{{\mathbb{N}}}

\newcommand{\As}{{\mathscr{A}}}\newcommand{\Bs}{{\mathscr{B}}}\newcommand{\Cs}{{\mathscr{C}}}
\newcommand{\Ds}{{\mathscr{D}}}
\newcommand{\Es}{{\mathscr{E}}}\newcommand{\Fs}{{\mathscr{F}}}
\newcommand{\Hs}{{\mathscr{H}}}
 
\newcommand{\Ms}{{\mathscr{M}}}\newcommand{\Os}{{\mathscr{O}}}
\newcommand{\Ps}{{\mathscr{P}}}	
\newcommand{\Rs}{{\mathscr{R}}}

\newcommand{\Xs}{{\mathscr{X}}}

\DeclareFontFamily{U}{rsfs}{\skewchar\font127 }
\DeclareFontShape{U}{rsfs}{m}{n}{%
   <5> <6> rsfs5
   <7> rsfs7
   <8> <9> <10> <10.95> <12> <14.4> <17.28> <20.74> <24.88> rsfs10
}{}
\DeclareSymbolFont{rsfs}{U}{rsfs}{m}{n} 
\DeclareSymbolFontAlphabet{\scr}{rsfs}

\newcommand{\Mf}{\scr{M}}

\DeclareMathOperator{\id}{Id} 

\DeclareMathOperator{\Ob}{Ob}

\DeclareMathOperator{\Hom}{Hom}
 
\DeclareMathOperator{\dom}{Dom}

\renewcommand{\emph}{\textbf} 										
\newcommand{\cj}[1]{\overline{#1}}									
\newcommand{\ip}[2]{\langle #1\mid #2\rangle}					
\renewcommand{\iff}{\Leftrightarrow}								

\newcommand{\hlink}[2]{\href{#1}{\texttt{#2}}} 

\newcommand{\xqedhere}[2]{%
  \rlap{\hbox to#1{\hfil\llap{\ensuremath{#2}}}}}

\theoremstyle{plain}

\theoremstyle{definition}

\numberwithin{equation}{section}  		
\title{\textbf{Categorical Non-commutative Geometry}
}

\author{\normalsize 
Paolo Bertozzini$^a$, 
Roberto Conti$^b$\footnote{Current Address: Dipartimento di Scienze di Base e Applicate per l'Ingegneria, Sezione di Matematica, 
Sapienza Universit\`a di Roma, Via A. Scarpa 16, I-00161 Roma, Italy. E-mail: \texttt{roberto.conti@sbai.uniroma1.it}} , 
Wicharn Lewkeeratiyutkul$^c$ 
\\
\normalsize $^a$\textit{Department of Mathematics and Statistics, Faculty of Science and Technology,}
\\
\normalsize \textit{Thammasat University, Bangkok 12121, Thailand} 
\\ 
\normalsize e-mail: \texttt{paolo.th@gmail.com}
\\ 
\normalsize $^b$\textit{Dipartimento di Scienze, Universit\`a di Chieti-Pescara ``G. D'Annunzio'',}
\\
\normalsize \textit{ Viale Pindaro 42, I-65127 Pescara, Italy} 
\\  
\normalsize e-mail: \texttt{conti@sci.unich.it}
\\
\normalsize $^c$\textit{Department of Mathematics, Faculty of Science, Chulalongkorn University,}
\\
\normalsize \textit{Bangkok 10330, Thailand}
\\ 
\normalsize e-mail: \texttt{Wicharn.L@chula.ac.th}
}

\date{\normalsize{08 December 2013}\footnote{This is a reformatted version for arXiv of an original paper submitted to the proceedings of the ``Algebra Geometry and Mathematical Physics 2010'' conference and published in: \textit{Journal of Physics: Conference Series} 346 (2012) 012003.}.}

\begin{document}

\maketitle

\begin{abstract} \noindent 
The purpose of this short note is to outline the current status of some recent research programs aiming at a categorification of parts of A.Connes 
non-commutative geometry and to provide an outlook on some possible future developments in categorical non-commutative geometry. 

\smallskip


\smallskip

\noindent
\emph{MSC-2010:} 						
					46L87,			
					18D05,			
					46M15, 			
					46M99,			
 					16D90.			

\end{abstract}

\tableofcontents




%
%
%
%
%
%

\section{Introduction}
Category theory has been applied to operator algebraic settings since many years ago, probably starting around 1980 with the pioneering work of John Roberts in algebraic quantum field theory~\cite{GLR} and since then it has been constantly used in the theory of superselection sectors. 

As far as we know, apart from a couple of independent proposals by A.Sitarz~\cite[section~3.2]{S} and Y.Manin~\cite{Ma}, the study of ``categorical non-commutative geometry'', in the setting of A.Connes' spectral triples, started around 2002-2003 as a by-product of our research project ``Modular Spectral Triples in Non-commutative Geometry and Physics'' (Thai Research Fund Grant: RSA4580030). There, in order to identify a non-commutative configuration space from  a non-commutative phase-space obtained by Tomita-Takesaki modular theory, a kind of polarization was necessary and for this purpose we were motivated to introduce a definition of sub-object and study the most elementary notion of morphism of spectral triples~\cite{BCL1}. 

The research on categorical non-commutative geometry soon started to become one of the main areas of our activity as documented in the survey paper~\cite{BCL2} that can still be considered a fairly good introduction to the subject. 

Apart from the study of alternative simple notions of morphism of spectral triples~\cite{BCL5}, we started a project of ``categorification'' of A.Connes non-commutative geometry. Around 2006 we introduced the terms ``horizontal categorification'' and ``vertical categorification'' in order to distinguish the categorical ``many-object version'' of usual mathematical concepts from the more demanding ``higher-morphism'' counterpart and we concentrated in proving a horizontal categorified version of Gel'fand-Na\u\i mark duality for commutative full C*-categories, a result that has recently appeared in~\cite{BCL3}. The spectrum of such a C*-category consists of a specific Fell line-bundle that we call ``topological spaceoid''. 

The horizontal categorification of Krein-C*-algebras (essentially categories of bounded linear operators between complete semi-definite linear spaces) has also been investigated in~\cite{BK}. 

In the meantime categorical non-commutative geometry (in A.Connes' sense) has been the subject of more and more investigations at different levels of technical sophistication by different authors:
\begin{itemize}
\item          
A.Connes, C.Consani, M.Marcolli in the fist part of their paper~\cite{CCM} describe two possible general approaches for the definition of categories of spectral triples, suggestions that are carried further by A.Connes, M.Marcolli, in their recent 
book~\cite{CM}, 
\item  
B.Mesland~\cite{M} introduced what is in our opinion the most sophisticated notion of category of spectral triples, based on ``smooth correspondences'' (certain specific KK-bimodules),  
\item 
R.Dawe Martins in~\cite{DM1,DM2,DM3} proposed other generalizations of the notion of spectral triple that are based on Fell bundles. 
\end{itemize}

Research in categorical non-commutative geometry (as already claimed in our original motivation) seems to be of particular interest in all the attempts to provide a formulation of algebraic quantum gravity (see~\cite{BCL4} and also~\cite{CM, DMZ, MZ, DM3}). 
 
\section{Morphisms of Spectral Triples}

Recall\footnote{For basic background on A.Connes' definition of spectral triples we refer to~\cite{C,GVF,V,CPR}.} that a (naive) spectral triple $(\As,\H,D)$ is given by a pre-C*-algebra $\As$ faithfully represented via bounded operators on a Hilbert space $\H$ and a possibly unbounded self-adjoint ``Dirac'' operator $D$ on $\H$ that has compact resolvent and  commutators $[D,\pi(x)]_-$, for all $x\in \As$, bounded on a common dense domain in $\H$.  
Typical examples of spectral triples originating in differential geometry are the Atiyah-Singer spectral triples of a compact orientable Riemannian spinorial manifold $M$, where $\As:=C^\infty(M)$ is the algebra of smooth complex-valued functions on $M$ represented by left multiplication on the Hilbert space $\H:=L^2(S(M))$ of square integrable section of a spinor bundle on $M$ with $D$ given by the usual Pauli-Atiyah-Singer Dirac operator. Compact orientable Riemannian manifolds can also be described by spectral triples taking $\As:=C^\infty(M)$ represented by left multiplication on the space $\H:=L^2(\Lambda_\bullet(M))$ of sections of the Grassmann bundle of $M$ with Dirac operator $D:=d+d^*$. 

\subsection{Totally Geodesic Morphisms}

In our first paper on this subject~\cite{BCL1} we proposed this notion of morphism: 
given two spectral triples $(\As_j,\H_j,D_j),$ with $j=1,2,$ a 
\emph{morphism of spectral triples} is a pair $(\phi,\Phi)$, 
where \hbox{$\phi: \As_1\to \As_2$} is a $*$-homomorphism between the pre-C*-algebras $\As_1,\As_2$ and $\Phi: \H_1\to \H_2$ is a bounded linear map that ``intertwines'' the representations $\pi_1, \pi_2\circ\phi$ and the Dirac operators $D_1, D_2$ i.e.: 
$\pi_2(\phi(x))\circ \Phi=\Phi\circ \pi_1(x)$, $\forall x \in \As_1$, 
$D_2\circ\Phi(\xi)= \Phi \circ D_1(\xi)$, $\forall \xi\in \dom D_1$. 

This definition of morphism clearly implies a strong relationship between the spectra of the Dirac operators of the two spectral triples. Loosely speaking, for $\phi$ epi and $\Phi$ coisometric (respectively mono and isometric), in the case of Atiyah-Singer spectral triples, one should expect such definition to become relevant only for maps that ``preserve the geodesic structures'' (totally geodesic immersions and respectively totally geodesic submersions).\footnote{Bertozzini P, Conti R, Lewkeeratiyutkul W, Non-commutative Totally Geodesic Submanifolds and Quotient Manifolds, work in progress.}  
Furthermore these morphisms depend, at least in some sense, on the spin structures:  
this ``spinorial rigidity'' (at least in the case of morphisms of real even spectral triples, when we also impose intertwining conditions between $\Phi$ and the real structures $J_j$ and the gradings $\Gamma_j$) requires that such  morphisms between spectral triples of different dimensions might be possible only when the difference in dimension is a multiple of~8. 

\subsection{Metric Morphisms} 
A notion of morphism that is essentially blind to the spin structures has been proposed in~\cite{BCL3} where it has been used to prove a refined version of Gel'fand duality for Atiyah-Singer spectral triples and metric isometries of spinorial manifolds.  
Given two spectral triples $(\As_j,\H_j,D_j)$, with $j=1,2$, denote by 
$d_{D_j}(\omega_1,\omega_2):=\sup\{|\omega_1(x)-\omega_2(x)| \ | \ x \in \As, \ \|[D_j,\pi(x)]_-\|\leq 1\}$
the quasi-distance induced on the sets $\Ps(\As_j)$ of pure states of $\As_j$.
A \emph{metric morphism} of spectral triples is a unital epimorphism\footnote{Note that if $\phi$ is an epimorphism, its pull-back 
$\phi^\bullet$ maps pure states into pure states.} 
$\phi : \As_1\to\As_2$ of pre-C*-algebras whose pull-back 
$\phi^\bullet : \Ps(\As_2)\to\Ps(\As_1)$, $\phi^\bullet(\omega):=\omega\circ \phi$ is an isometry, i.e.
$d_{D_1}(\phi^\bullet(\omega_1),\phi^\bullet(\omega_2))=d_{D_2}(\omega_1,\omega_2)$, for all 
$\omega_1,\omega_2\in \Ps(\As_2)$.

\subsection{Riemannian Morphisms} 

A weaker notion of metric morphisms (that in the case of isomorphisms reduces to the unitary maps considered in~\cite{PV}) and that for Atiyah-Singer spectral triples should reproduce the usual situation of Riemannian immersions and submersions of spinorial manifolds is as follows:\footnote{Bertozzini P, Conti R and Lewkeeratiyutkul W, Categories of Spectral Triples and Morita Equivalence, work in progress.} given two spectral triples $(\As_j,\H_j,D_j),$ with $j=1,2,$ a 
\emph{Riemannian morphism} is a pair $(\phi,\Phi)$
where $\phi: \As_1\to \As_2$ is a $*$-homomorphism between the pre-C*-algebras $\As_1,\As_2$ and $\Phi: \H_1\to \H_2$ is a bounded linear map that ``intertwines'' the representations $\pi_1, \pi_2\circ\phi$ and the commutators of the Dirac operators 
$D_1, D_2$: 
$\pi_2(\phi(x))\circ \Phi=\Phi\circ \pi_1(x)$, $\forall x \in \As_1$, 
$[D_2,\pi_2(\phi(x))]_-\circ\Phi= \Phi \circ [D_1,\pi_1(x)]_-$, $\forall x\in \As_1$. 
Note that the boundedness of $\Phi$ (here as well as in the case of totally geodesic morphisms) can actually be weakened, considering unbounded operators, the important property here being the fact that the adjoint action of $\Phi$ on the algebra 
$\Omega_D(\As)$ generated by $\pi(\As)$ and the commutators $[D,\pi(x)]_-$, $x\in \As$ (the ``non-commutative Clifford algebra'') is a $*$-homomorphism extending $\phi$.

\subsection{Morita Morphisms of Spectral Triples}

All of the several definitions of morphisms considered above have been essentially modelled on the case of commutative algebras of functions, where $*$-homomorphisms are abundant, and although they still make sense in the non-commutative case, they correspond to quite special ``maps'' of non-commutative spaces. 
In a wider perspective, a morphism of spectral triples $(\As_j,\H_j,D_j)$, for $j=1,2$, should be formalized as a ``suitable'' functor 
$\Fs: {}_{\As_1}\Mf\to {}_{\As_2}\Mf$, between the categories ${}_{\As_j}\Mf$ of $\As_j$-modules, having ``appropriate intertwining'' properties with the Dirac operators $D_j$. 
Under some ``mild'' hypothesis, by Eilenberg-Gabriel-Watt theorem, any such functor is given by ``tensorization'' with a bimodule. These bimodules, suitably equipped with spectral data (as in the case of spectral triples), will provide the natural setting for a general theory of morphisms of non-commutative spaces. This ``Morita morphism'' point of view has been first advocated by Y.Manin~\cite{Ma}, but it is had already been implicitly exploited in A.Connes' ``transfer'' of Dirac operators via Morita equivalence bimodules equipped with a connection~\cite{C2,CC}.  

In~\cite{BCL2} we also noticed the construction of a strictly related category of \emph{Morita-Connes morphisms} of spectral triples (containing A.Connes' ``transfers and inner deformations'' as isomorphisms) based on the choice of a connection on a Morita morphism (that is not necessarily an imprimitivity bimodule) i.e.: a left-$\As_2$ right-$\As_1$ bimodule that is a Hilbert C*-module over $\As_1$, a Hermitian connection\footnote{Here $\Omega^1_D(\As)$ denotes the $\As$-bimodule inside the algebra $\Omega_D(\As)$ spanned by the commutators $[D,\pi(x)]_-$, $x\in \As$.} 
$\nabla:X\to X\otimes_{\As_1}\Omega^1_{D_1}(\As_1)$ on the bimodule $X$ 
(the Dirac operators on the spectral triples $(\As_j,\H_j,D_j)$, $j=1,2$, being related to the connection $\nabla$ by the Connes' ``transfer'' formula $D_2(\xi\otimes h)=\xi\otimes D_1(h)+(\nabla\xi)(h)$ where $h\in \H_1$ and $\xi\in X$) 
and with composition given by the bimodule $X^3:=X^2\otimes_{\As_2} X^1$ equipped with the connection:
$\nabla^3(\xi_2\otimes\xi_1)(h):=\xi_2\otimes(\nabla^1\xi_1)(h)+(\nabla^2\xi_2)(\xi_1\otimes h)$, where $\xi_1\in X^1,\xi_2\in X^2$, 
$h\in \H_1$. 

\subsection{Mesland Morphisms}

Morphism of spectral triples via Morita correspondences have been further developed in the works by 
A.Connes, M.Marcolli~\cite[chapter~8.4]{CM} and M.Marcolli, A.Z.al Yasri~\cite{MZ} were ``spectral correspondences'', defined as Hilbert C*-bimodules, are used to provide a ``bivariant version'' of spectral triples. 

The most complete proposal in this direction comes from the work by B.Mesland~\cite{M} that has defined a category of (unitary equivalence classes of) smooth KK-bimodules that seems to be the best candidate for a non-commutative metric category of spectral triples. 
A \emph{Mesland morphisms} from the  spectral triple $(\Bs,\Hs',D')$ to the spectral triple $(\As,\Hs,D)$ is given by a unitary isomorphism class of an unbounded ``smooth'' $\As$-$\Bs$-bimodule $(\Es,S)$ with ``smooth connection'' $\nabla$ such that:
$\Hs$ is isomorphic to $\Es\cj{\otimes}_\Bs\Hs'$ (where here $\cj{\otimes}$ denotes the Haagerup tensor product);  
$[\nabla,S]$ is a completely bounded operator; 
$D=S\otimes \id + \id\otimes_\nabla D'$ with $\id\otimes_\nabla D'(x\otimes h):=(-1)^{\partial x}(x\otimes D'h+(\nabla_{D'}x)h)$, where $x\in \Es$, $h\in \H'$ ($\partial x$ denoting the degree of $x$ in the graded module $\Es$). 

\section{Categorification of Gel'fand Na\u\i mark Duality}

\subsection{Horizontal Categorification} 

In the same way as a category can be seen as a ``many-objects'' version of a monoid or a groupoid can be thought as a multi-objects version of a group, a C*-category is a (horizontal) categorification of a C*-algebra. Furthermore, in the same way as every category $\Cs$ induces a projection functor $\pi:\Hom_\Cs\to\Ob_\Cs\times\Ob_\Cs$, a C*-category $\Cs$ can be identified as a very special kind of Fell bundle where the base category is a ``double groupoid'' $\Ob_\Cs\times \Ob_\Cs$.  

More precisely, given an \emph{inverse involutive category} $\Xs$ (i.e.~a category equipped with an object-preserving contravariant functor $x\mapsto x^*$ such that for all arrows $x\in \Xs$, $(x^*)^*=x$ and such that $x\circ x^*\circ x=x$ for all $x\in \Xs$) a \emph{unital Fell bundle over $\Xs$} is a Banach bundle\footnote{We refer to J.Fell, R.Doran~\cite{FD} for all the details on Banach bundles.} $\pi:\Es\to\Xs$ with, a total space $\Es$ that is 
an involutive category, a projection $\pi$ that is a covariant $*$-functor and such that: the composition in $\Es$ is fiberwise bilinear and norm submultiplicative; the involution in $\Es$ is fiberwise conjugate linear with the C*-property $\|e^*\circ e\|=\|e\|^2$ and such that $e^*\circ e$ is a positive element in the C*-algebra $\pi^{-1}(\pi(e^*\circ e))$.\footnote{Note that $\pi^{-1}(x^*\circ x)$ is always a unital C*-algebra and $\pi^{-1}(x)$ is always a C*-bimodule onto the C*-algebras $\pi^{-1}(x^*\circ x)$ and 
$\pi^{-1}(x\circ x^*)$.}  
A Fell bundle is \emph{saturated} whenever the Hilbert bimodules 
$\pi^{-1}(x)$ are full over the C*-algebras $\pi^{-1}(x^*\circ x)$ and $\pi^{-1}(x\circ x^*)$. 
A (small) \emph{C*-category} can be identified as a unital Fell bundle over an involutive category of the form $\Os\times\Os$ for a certain set $\Os$. It is said to be full if it is saturated as a Fell bundle and commutative if the C*-algebras $\pi^{-1}(x^*\circ x)$ are Abelian.  

In the search for an appropriate notion of ``spectrum of a commutative full small C*-category'' we defined a \emph{topological spaceoid} as a unital Fell bundle of rank-one (i.e.~with one-dimensional fibers) whose base category is given by a direct product $\Delta_X\times \Rs_\Os$, where $X$ is a compact Hausdorff space, $\Os$ is a discrete space, $\Delta_X:=\{(x,x) \ | \ x\in X\}$ is the ``diagonal of $X$'' and $\Rs_\Os:=\Os\times\Os$ is the maximal equivalence relation on $\Os$.  

In~\cite{BCL3}, we provided a categorical extension of the usual Gel'fand-Na\u\i mark duality between the category of unital $*$-homomorphisms of Abelian unital C*-algebras and the category of continuous maps of compact Hausdorff spaces to a new duality between the category of object-preserving \hbox{$*$-functors} of small commutative full C*-categories and a category of suitable morphisms of space\-oids. 

In~\cite{BCL6} we further generalized the notion of Fell bundle introducing a definition of involutive categorical bundle  (Fell bundle) enriched in an involutive monoidal category (or even in an involutive 2-fold category) and we made use of this concept to relate three equivalent ways to describe the spectrum of a full commutative small C*-category. 

An interesting by-product of this investigation is an alternative direct proof of a spectral theorem for imprimitivity Hilbert C*-bimodules over Abelian C*-algebras~\cite{BCL7} (i.e.~a Hermitian version of Serre-Swan theorem) that is suitable to provide a ``bivariant version'' of A.Takahashi's duality between categories of Hilbert C*-modules and categories of Hilbert 
bundles~\cite{Ta1,Ta2}. 

\subsection{Non-full C*-categories}

One further essential step is to extend our spectral theorem to the case of non-full small commutative 
C*-categories. In this case the spectrum of the C*-category $\Cs$ is no more a line-bundle and can be described as a Fell bundle with fibers of dimension less than or equal to one. The locus of base points supporting zero-dimensional fibers is given by a family of closed sets $F_{AB}$ for all $A,B\in \Ob_\Cs$ with the properties $F_{AA}=\varnothing$, $F_{AB}=F_{BA}^*$ and 
$F_{AC}\subset F_{AB}\circ F_{BC}$ for all $A,B,C\in \Ob_\Cs$. These non-full categories correspond of course to special cases of closed ``ideals'' of full commutative C*-categories. 

\subsection{Vertical Categorification}

In view of a further vertically categorified extension of Gel'fand duality, 
we are investigating the existence of reasonable notions 
of strict $n$-C*-categories.\footnote{Bertozzini P, Conti R, Lewkeeratiyutkul W, Suthichitranont N, Strict Higher C*-categories, work in progress. 
\\ 
See also the slides ``Categories of Non-commutative Geometries'' at the second workshop 
``\href{http://categorieslogicphysics.wikidot.com/meeting2}{Categories, Logic and Physics}'' in Imperial College.}

Recall~\cite[section~1.4]{L} that a \emph{globular $n$-set}
$\Cs^0\leftleftarrows\Cs^1\leftleftarrows \cdots \Cs^{m-1}\leftleftarrows\Cs^m\leftleftarrows\dots \leftleftarrows \Cs^n$, 
$n\in \NN$, is given by:
a collections of classes $\Cs^m$, for all $m=0,\dots,n$, whose elements are called \emph{$m$-arrows}, 
and a pair of \emph{source, target} maps 
$s_m,t_m:\Cs^m\to \Cs^{m-1}$, for all $m=1,\dots,n$, such that for all $m=1,\dots, n-1$, we have 
$s_m\circ s_{m+1}=s_m\circ t_{m+1}$,   
and $t_m\circ s_{m+1}=t_m\circ t_{m+1}$.

A (globular) \emph{strict $n$-category} (for example see T.Leinster~\cite[section~1.4]{L}) has been defined as a globular $n$-set that for all $0\leq p<m\leq n$, is equipped with a partial \emph{$p$-composition} map 
$\circ^m_p: \Cs^m \times_{\Cs^p} \Cs^m \to \Cs^m$, $(x,y)\mapsto x\circ^m_p y$, 
defined on the set $\Cs^m\times_{\Cs^p}\Cs^m$ of $p$-composable $m$-arrows 
$(x,y)\in\Cs^m\times_{\Cs^p}\Cs^m \iff t_{p+1}\circ\cdots\circ t_m(y)=s_{p+1}\circ\cdots\circ s_m(x)$, 
such that, for all $m=0,\dots, n-1$, there is an \emph{identity map}   
$\iota_m:\Cs^m\to\Cs^{m+1}$, 
in such a way that the following axioms are satisfied:
\begin{itemize}
\item 
for all $m=0,\dots,n$, for all $p=0,\dots,m-1$, for all $(x,y)\in \Cs^m\times_{\Cs^p}\Cs^m$, \\ 
$s_m(x\circ_p^m y)=s_m(y), \quad t_m(x\circ_p^m y)=t_m(x), \quad  \text{if $p=m-1$}$;  \\
$s_m(x\circ_p^m y)=s_m(x)\circ^{m-1}_p s_m(y), \quad 				 \text{if $p=0,\dots,m-2$}$, \\ 
$t_m(x\circ_p^m y)=t_m(x)\circ^{m-1}_p t_m(y), \quad 				 \text{if $p=0,\dots,m-2$}$;
\item 
for all $x\in \Cs^m$, 
$s_{m+1}(\iota_m(x))=x, \quad t_{m+1}(\iota_m(x))=x$; 
\item 
for all $m=1,\dots,n$ and $p=0,\dots,m-1$ and for all $x,y,z\in \Cs^m$,  
$(x\circ_p^m y) \circ_p^m z = x\circ_p^m (y \circ_p^m z)$, 
whenever $(x,y),(y,z)\in \Cs^m\times_{\Cs^p}\Cs^m$ holds;  
\item 
for all $m=1,\dots,n$, for all $p=0,\dots,m-1$, for all $x\in \Cs^m$, 
\\ 
$\Big(\iota_{m-1} \circ \dots \circ \iota_p \big(t_{p+1}\circ \dots \circ t_m(x)\big)\Big) \circ^m_p x = x$, \quad
$x= x\circ^m_p \Big(\iota_{m-1} \circ \dots \circ \iota_p \big(s_{p+1}\circ \dots \circ s_m(x)\big)\Big)$;  
\item 
for all $m=2,\dots,n$, for all $p,q=0,\dots, m-1$, with $q<p$, for all $w,x,y,z\in \Cs^m$ such that 
$(w,x),(y,z)\in \Cs^m\times_{\Cs^p}\Cs^m$ and $(w,y),(x,z)\in \Cs^m\times_{\Cs^q}\Cs^m$, we have the exchange property 
$(w\circ^m_p x)\circ^m_q (y\circ^m_p z)=(w\circ^m_q y)\circ^m_p (x\circ^m_q z)$;
\item 
for all $m=1,\dots,n-1$, for all $p=0,\dots,m-1$, for all $(x,y)\in \Cs^m\times_{\Cs^p}\Cs^m$, we have 
$\iota_m(x\circ^m_p y)=\iota_m(x)\circ^{m+1}_p\iota_m(y)$.
\end{itemize}

It is reasonable to define a \emph{strict involutive $n$-category} as a strict $n$-category that is equipped with 
a family of ``involutions'' $*^m:\Cs^m\to\Cs^m$, for $0<m\leq n$, that satisfy the following properties:\footnote{
Actually it is perfectly possible to require the existence of involutions only for certain specific ``arrow levels'' so that, in the case of involution present only for the level $m=n$, previous definitions of $2$-C*-categories can be recovered. In the opposite direction, it might also be possible to require further axioms for involutions $*^m_q:\Cs^m\to\Cs^m$ of depth $q$ for $0\leq q<m\leq n$, but we will not go into further details here.
}  
\begin{itemize}
\item   
$s_m(x^{*^m})=t_m(x), \quad t_m(x^{*^m})=s_m(x)$, for all $x\in \Cs^m$,\item 
$(x\circ^m_p y)^{*^m}=y^{*^m}\circ^m_p x^{*^m}, \ \text{for $p=m-1$}$,  \quad 
$(x\circ^m_p y)^{*^m}=x^{*^m}\circ^m_p y^{*^m}, \ \text{for $0\leq p<m-1$}$, 
for all $x,y\in \Cs^m\times_{\Cs^p}\Cs^m$ with $m=1,\dots,n$, 
\item  $(x^{*^m})^{*^m}=x$, for all $x\in \Cs^m$.  
\end{itemize}

Finally, one might try to define a \emph{strict-$n$-C*-category} to be a strict involutive $n$-category such that: 
\begin{itemize}
\item
for all $m=1,\dots,n$, and $x,y\in \Cs^{m-1}$, the sets $\Cs^m(x,y):=\{z\in\Cs^m \ | \ s_m(z)=y, \ t_m(z)=x\}$ 
are Banach spaces with norm denoted by $x\mapsto\|x\|_m$, 
for $0\leq p<m$, 
\item 
for all $w,x,y,z\in \Cs^{m-1}$ such that 
$\Cs^m(w,x)\times\Cs^m(y,z)\subset \Cs^m\times_{\Cs^p}\Cs^m$, 
the composition maps $\circ^m_p:\Cs^m(w,x)\times \Cs^m(y,z)\to\Cs^m$ are bilinear, 
\item
for all $m=1,\dots,n$, for all $x,y\in\Cs^{m-1}$, the maps $*^m:\Cs^m(x,y)\to \Cs^m$ are conjugate linear; 
\item 
for all $m=1,\dots,n$, for all $p=0,\dots,m-1$, for all pairs $(x,y)\in \Cs^m\times_{\Cs^p}\Cs^m$,
\\ 
$\|x\circ_p^m y\|_m\leq \|x\|_m\cdot\|y\|_m$,
\item 
for all $m=1,\dots,n$ and $0\leq p<m$, for all $(x^{*^m},x)\in \Cs^m\times_{\Cs^p}\Cs^m$, 
$\|x^{*^m}\circ^m_p x\|_m=\|x\|_m^2$.
\end{itemize}
Note that the above properties already imply that, for all $m=1,\dots,n$ and for all $x\in \Cs_{m-1}$, the set $\Cs^m(x,x)$ is a C*-algebra with multiplication $\circ^m_{m-1}$ and involution $*^m$ and hence the following final condition is meaningful:
\begin{itemize} 
\item
for all $m=1,\dots,n$, for all $x\in \Cs^m(u,v)$, 
$x^{*^m}\circ^m_{m-1}x\in\Cs^m(u,u)_+$, i.e.~$x^{*^m}\circ^m_p x$ is a positive element in the C*-algebra $\Cs^m(u,u)$. 
\end{itemize}

A \emph{left module ${}_\Cs\Ms$ over the $n$-category $\Cs$} is given by
\begin{equation*}
\xymatrix{
\Cs^0 & \ar[l]^s_t \Cs^1 & \ar[l]^s_t \cdots & \ar[l]^s_t \Cs^{n-1} & \ar[l]^s_t \Cs^n 
\\
 & \Ms^1 \ar[ul]^\tau & \Ms^2 \ar[ul]^\tau & \cdots & \Ms^{n} \ar[ul]^\tau & 
}
\end{equation*} 
where for all $m=0,\dots,n$, $\tau:\Ms^m\to\Cs^{m-1}$ is a fibered category over the $(m-1)$-category $\Cs^{m-1}$ and, 
for all $0\leq p<m\leq n$, there is a left action $\mu^m_p:\Cs^m\times\Ms^m\to\Ms^m$ of the bi-fibered $(m-1)$-category 
$\Cs^m\rightrightarrows\Cs^{m-1}\times\Cs^{m-1}$ over $\Ms^m\to\Cs^{m-1}$ such that 
$\mu^m_p(\Cs^m(x,y)\times\Ms^m(z))\subset \Ms^m(x)$ whenever $(y,z)\in\Cs^{m-1}\times_{\Cs^p}\Cs^{m-1}$ with 
$x=y\circ^{m-1}_p z$.\footnote{For $p=m-1$ we assume $\Cs^{m-1}\times_{\Cs^p}\Cs^{m-1}=\Delta_{\Cs^{m-1}}$.} 
Similar definitions can be given for right modules $\Ms_\Cs$ and bimodules ${}_\Cs\Ms_\Cs$ over the $n$-category $\Cs$. 

The notion of \emph{left Hilbert C*-module ${}_\Cs\Ms$ over a strict $n$-C*-category $\Cs$} should be given imposing that for all $m=1,\dots n$, $\tau:\Ms^m\to\Cs^{m-1}$ is a ``Fell bundle''(for all the compositions and involutions in $\Cs^{m-1}$) equipped with an inner product $\ip{\cdot}{\cdot}_m:\Ms^m\times\Ms^m\to\Cs^m$ such that $\ip{\Ms^m(x)}{\Ms^m(y)}_m\subset\Cs^m(y,x)$.\footnote{Again, corresponding definitions can be given for right Hilbert C*-modules and right/left bimodules over a strict \hbox{$n$-C*-category}, but it will be necessary to distinguish right and left structures also for bimodules.}  

Examples of rank-one strict-$n$-C*-categories i.e.~strict-$n$-C*-categories such that the Banach space $\Cs^m(x,y)$ is one-dimensional, for every \hbox{$m=1,\dots,n$}, 
can be constructed by hand recursively. In the theory of higher C*-categories they play the role of the scalar field $\CC$.
Hilbert C*-modules over rank-one strict $n$-C*-categories will play the role of $n$-Hilbert spaces. 
Examples of non-commutative strict-$n$-C*-categories are expected to arise as 
 ``categories of endomorphisms'' of left Hilbert C*-modules over rank-one $n$-C*-categories.

A formulation of Gel'fand-Na\u\i mark duality in such higher C*-categorical context requires the usage of ``iterated Fell line-bundles'' and it is under investigation. 

\subsection{Horizontal Categorification of Spectral Triples}

One of the main original motivations in the study of C*-categories comes from the realization that, since the ``off-diagonal blocks'' 
$\Cs_{AB}:=\Hom_\Cs(B,A)$ are Hilbert C*-bimodules over the C*-algebras $\Cs_{AA}$ and $\Cs_{BB}$, the study of possible axiomatizations of spectral triples over \hbox{C*-categories} might provide some further light on the appropriate definition of ``bivariant spectral triples'' and more generally Morita morphisms of spectral triples. Of course in the case of full C*-categories, all the bimodules $\Cs_{AB}$ are imprimitivity bimodules (i.e.~isomorphisms in the Morita-Rieffel category of Hilbert C*-bimodules) and so, in this special case, we are bound to obtain arrows in a groupoid of isomorphisms of spectral triples. Without entering into further details that will be developed elsewhere, we note that spectral triples over a C*-category can be simply defined as spectral triples over the enveloping C*-algebra of the C*-category. 
For example, in the attempt to generalize naive spectral triples to a categorified context\footnote{Bertozzini P, Conti R and Lewkeeratiyutkul W, Spectral Geometries over C*-categories and Morphisms of Spectral Geometries, work in progress.} we can define a \emph{categorical spectral geometry} as a triple $(\Cs,\H,\Ds)$ given by: 
\begin{itemize}
\item
a pre-C*-category $\Cs$; 
\item
a module $\Hs$ over $\Cs$ that is also a Hilbert C*-module over $\CC$; in other terms a family of Hilbert spaces $\Hs$ equipped with an object bijective $*$-functor $\pi:\Cs\to\Bs(\Hs)$ with values in the C*-category of bounded linear maps between the Hilbert spaces in the family $\Hs$; 
\item
the generator $\Ds$ of a unitary one-parameter group on $\Hs$ (i.e.~the generator of a one-parameter group whose adjoint action in the enveloping C*-algebra of $\Bs(\Hs)$ leaves $\Bs(\Hs)$ invariant) such that, for all $x\in \Cs$, 
$[\Ds,\pi(x)]_- $ is extendable to an operator in $\Bs(\Hs)$. 
\end{itemize} 

In the case of C*-categories, the notion of bimodule over a C*-category is significantly different from that of left or right module (see for example P.Mitchener~\cite{Mi}) and this results in a further complication as can be seen in the following very tentative definition. 

A \emph{bivariant spectral geometry} over two pre-C*-categories (with the same objects) $\As$ and $\Bs$ is a quintuple 
$(\As,\Bs,\Hs,\Ds_\As,\Ds_\Bs)$, where 
\begin{itemize}
\item 
$\Hs$ is a bimodule over $\As$-$\Bs$ that is also a Hilbert C*-bimodule over $\CC$ and hence it is equipped with two 
$*$-representations $\rho:\As\to\Bs_\rho(\Hs)$ and $\lambda:\Bs\to\Bs_\lambda(\Hs)$ into the right, and respectively the left, 
C*-category of the bimodule; 
\item 
$\Ds_\As$ (acting on the left) and $\Ds_\Bs$ (acting on the right) are two (generally unbounded) self-adjoint operators on $\Hs$ that generate on the enveloping algebras of $\Bs_\rho(\Hs)$, and respectively of $\Bs_\lambda(\Hs)$, one-parameter groups leaving $\Bs_\rho(\Hs)$, and respectively $\Bs_\lambda(\Hs)$, invariant and such that $[\Ds_\As,\rho(x)]_-$ and $[\Ds_\Bs,\lambda(y)]_-$ are extensible to bounded operators in $\Bs_\rho(\Hs)$, $\Bs_\lambda(\Hs)$, for all $x\in \As$ and $y\in \Bs$. 
\end{itemize}

\section{Outlook}

A short-term objective of this line of research is to provide explicit examples of functors from suitable categories of geometrical spaces (such as for example oriented Riemannian or spinorial compact manifolds) to categories of spectral triples (such as the category described by B.Mesland~\cite{M} or possibly some variants of it). 

Since spectral triples are a very sophisticated kind of mathematical tool where topological, measurable, smooth and metric structures are simultaneously present, it seems worth to spend some time investigating separately the categorical structures involved in the case of oriented spaces, measure spaces, (Riemannian/Hermitian) manifolds/bundles equipped with connections or with spinorial bundles and their ``dual'' categories of modules. 
Of particular interest is the case of ``non-commutative measure spaces'' and the study of the categorical structure implicit in Tomita-Takesaki modular theory and in Falcone-Takesaki non-commutative flow of weights.  

Some more ambitious goals include: 
\begin{itemize}
\item
spectral reconstruction theorems for certain classes of morphism of spectral triples,  
\item 
extensions of our Gel'fand duality result to full non-commutative C*-categories, 
\item 
a ``spectral theory'' of spectral triples in terms of Fell line-bundles (along the lines envisaged by R.Martins) and its application to physics, 
\item 
the study of possible relations between (categorical) non-commutative geometry and \\ Grothendieck's topoi.  
\end{itemize}
Some applications of such mathematical structures in physics are also priorities: 
\begin{itemize}
\item 
in the context of loop quantum gravity, we might provide (categorical) non-commutative geometries associated to the ``quantum geometries'' described by spin-networks,
\item 
usage of (higher/modular) categorical structures to obtain a mathematical formulation of C.Rovelli's relational quantum mechanics, 
\item 
further progress in our modular algebraic quantum gravity proposal (see~\cite{BCL4}). 
\end{itemize}

\end{document}